\documentclass[12pt,a4paper]{article}

\usepackage{amsfonts}
\usepackage{amsmath}
\usepackage{amssymb}
\usepackage{amsthm}
\usepackage{array,booktabs}
\usepackage{color}
\usepackage[english]{babel}
\usepackage{epic}
\usepackage{epsfig}
\usepackage{enumitem}
\usepackage{graphics}
\usepackage{graphicx}
\usepackage{hyperref}
\usepackage{latexsym}
\usepackage{multicol}
\usepackage{multicol}
\usepackage{pstricks}
\usepackage{soul}
\usepackage{tikz}
\tikzset{
	every node/.style={circle, inner sep=2pt}
}
\usetikzlibrary{matrix,arrows,calc}
\usepackage[utf8]{inputenc}
\usepackage{url}

\newtheorem{theorem}{Theorem}[section]
\newtheorem{proposition}[theorem]{Proposition}
\newtheorem{lemma}[theorem]{Lemma}
\newtheorem{corollary}[theorem]{Corollary}
\newtheorem{definition}[theorem]{Definition}

\newtheorem{conjecture}{Conjecture}

\newcommand{\Z}{\mathbb{Z}}

\def\vec0{\mbox{\boldmath $0$}}

\def\D{\mbox{\boldmath $D$}}

\def\dist{\mbox{\rm dist}}

\begin{document}
	
	\title{On bipartite biregular large graphs\\
		derived from difference sets\footnote{
			The research of the G. Araujo-Pardo is supported by PAPIIT-M{\'e}xico under Project IN101821 and Proyecto de Intercambio Acad\'emico 1926: ``Gr\'aficas de Moore, Jaulas, Diseños de Bloques y N\'umero Crom\'atico'', CIC-UNAM. The research of C. Dalf\'o, M. A. Fiol, and N. L\'opez has been supported by
			AGAUR from the Catalan Government under project 2021SGR00434 and MICINN from the Spanish Government under project PID2020-115442RB-I00.
			The research of M. A. Fiol was also supported by a grant from the  Universitat Polit\`ecnica de Catalunya with references AGRUPS-2022 and AGRUPS-2023.}}
	\author{G. Araujo-Pardo$^a$, C. Dalf\'o$^b$, M. A. Fiol$^c$, N. L\'opez$^d$\\
		\\
		\small{$^a$Instituto de Matem\'aticas, Universidad Nacional Aut\'onoma de M\'exico,}\\
		{\small Mexico, \texttt{garaujo@math.unam.mx}}\\
		\small{$^b$Departament de Matem\`{a}tica, Universitat de Lleida,}\\
		\small{Igualada (Barcelona), Catalonia, \texttt{cristina.dalfo@udl.cat}}\\
		\small{$^c$Departament de Matem\`{a}tiques, Universitat Polit\`{e}cnica de Catalunya,}\\
		\small{Barcelona Graduate School of Mathematics,}\\
		\small{Institut de Matem\`atiques de la UPC-BarcelonaTech (IMTech)}\\
		\small{Barcelona, Catalonia, \texttt{miguel.angel.fiol@upc.edu}}\\
		\small{$^d$ Departament de Matem\`atica, Universitat de Lleida,}\\
		\small{Lleida, Spain, \texttt{nacho.lopez@udl.cat}}}
	
	\date{}
	
	\maketitle
	
	%\newpage
	%\tableofcontents

	%%%%%%%%%%%%%%%%%%%%%%%%%%%%%%
	%----------------Preliminaries---------------------------------
	%%%%%%%%%%%%%%%%%%%%%%%%%%%%%%
	\begin{abstract}
		
		A bipartite graph $G=(V,E)$ with $V=V_1\cup V_2$ is biregular if all the vertices of each stable set, $V_1$  and $V_2$, have the same degree, $r$ and $s$, respectively. This paper studies difference sets derived from both Abelian and non-Abelian groups. From them, we propose some constructions of bipartite biregular graphs with diameter $d=3$ and asymptotically optimal order for given degrees $r$ and $s$.
		Moreover, we find some biMoore graphs, that is, bipartite biregular graphs that attain the Moore bound.
		%a large number of vertices $N(r_1,r_2;d)$, together with their spectra.
		
	\end{abstract}

	%%%%%%%%%%%%%%%%%%%%%%%%%%%%%%%%%%%%%%%%%%%%%%%%%%%%%%%%%%%%%%%%%%%
	\section{Introduction}
	
	The \emph{degree/diameter problem} for graphs is finding the largest order of a graph with the prescribed degree and diameter. The maximum of this number is the \emph{Moore bound}, and a graph whose order coincides with this bound is called a {\emph{Moore graph}}.
	There is much work related to this topic (see the survey by Miller and \v{S}ir\'a\v{n} \cite{ms16}), and also to this subject with some restrictions of the original problem. One of them is related to the bipartite Moore graphs. In this case, the goal is to find regular bipartite graphs with maximum order and fixed diameter. In this paper, we deal with the problem initiated by Yebra, Fiol, and  F\`abrega \cite{yff83} in 1983, which consists of finding biregular bipartite Moore graphs.
	
	A bipartite graph $G=(V,E)$ with $V=V_1\cup V_2$ is {\em biregular} if all the vertices of the stable set $V_i$, for $i=1,2$, have the same degree. We denote $[r,s;d]$-bigraph a bipartite biregular graph of degrees $r$ and $s$ and diameter $d$; and by  $[r,s;d]$-biMoore graph the bipartite biregular graph of diameter $d$ that attains the Moore bound, which is denoted    $M(r,s;d)$.
	Notice that constructing these graphs is equivalent to constructing block designs, where one partite set corresponds to the points of the block design, and the other set corresponds to the blocks of the design. Moreover, each point is in a fixed number $s$ of blocks, and each block's size equals $r$. The incidence graph of this block design is an $[r,s;d]$-bigraph.
	
	We propose some constructions of bipartite biregular graphs with diameter $d=3$ and asymptotically optimal order for given degrees $r$ and $s$. Our main tool is the use of difference sets. The same approach was used by Erskine, Fratri\v{c}, and \v{S}ir\'a\v{n} \cite{efs21}, where they employed perfect difference sets to obtain graphs with diameter two and asymptotically optimal order for their maximum degree. Moreover, they proved that their graphs were isomorphic to some of the Brown graphs \cite{b66}, a well-known family of graphs in the degree/diameter problem.
	
	This paper is structured as follows. The following subsection gives the Moore-like bound for bipartite biregular graphs with an odd diameter. In Subsection \ref{sec:Moore-millorada-d=3}, we give a new improved Moore bound for diameter 3. In Section \ref{sec:d=3}, we study the bipartite biregular graphs derived from difference sets of Abelian groups. We give the so-called covering difference sets for $s=7$ in Section \ref{sec:d=4}. Finally, in Subsection \ref{sec:difference-sets-non-Abelian}, we deal with a bipartite biregular graph derived from difference sets of non-Abelian groups.
	
	\subsection{Moore-like bounds}
	\label{sec:Moore-like}
	Let $G=(V,E)$, with $V=V_1\cup V_2$, be an $[r,s;d]$-bigraph, where each vertex of $V_1$ has degree $r$, and each vertex of $V_2$ has degree $s$. Let $N_i=|V_i|$ for $i=1,2$. Then, counting in two ways the number of edges of $G$, we have
	$r N_1=s N_2$.
	Moreover, if the diameter is odd, say, $d=2m+1$ (for $m\ge 1$), and $u\in V_1$, we have (from the number of vertices at distance $0,1,2,\ldots,d-1$ from $u$):
	\begin{equation}
		\label{N1-odd}
		N_1\le 1+r(s-1)\frac{[(r-1)(s-1)]^{m}-1}{(r-1)(s-1)-1}=N_1',
	\end{equation}
	whereas, if $u\in V_2$ (interchanging $r$ and $s$):
	\begin{equation}
		\label{N2-odd}
		N_2\le 1+s(r-1)\frac{[(r-1)(s-1)]^{m}-1}{(r-1)(s-1)-1}=N_2'.
	\end{equation}
	Since $N_1'r\neq N_2's$, the Moore bound must be smaller than $N_1'+N_2'$. It was proved in Yebra, Fiol, and F\`abrega \cite{yff83} that, assuming $r>s$,
	\begin{equation}
		\label{N1-N2-odd}
		N_1\le \left\lfloor \frac{N_2'}{\rho}\right\rfloor \sigma\qquad\mbox{and}
		\qquad N_2\le \left\lfloor \frac{N_2'}{\rho}\right\rfloor \rho,
	\end{equation}
	where $\rho=\frac{r}{\gcd\{r,s\}}$ and $\sigma=\frac{s}{\gcd\{r,s\}}$.
	Thus, the Moore bound for odd diameter $d=2m+1$ is
	\begin{equation}
		\label{Moore-odd}
		M(r,s;2m+1)= \left\lfloor\frac{1+s(r-1)\frac{[(r-1)(s-1)]^{m}-1}{(r-1)(s-1)-1}}{\rho}\right\rfloor (\rho+\sigma).
	\end{equation}

	%See Tables \ref{tab:d=3} and \ref{tab:d=5} for  $2\le s\le r\le 10$ and $d=3$ and $d=5$, respectively. The known attainable bounds are in boldface.
	
	For diameter $d=3$ ($m=1$), the Moore bounds in \eqref{Moore-odd} become 
	\begin{equation}
		M(r,s;3)=\left\lfloor \frac{1+s(r-1)}{\rho} \right\rfloor(\rho+\sigma),
		\label{Moore-odd(d=3)}
	\end{equation}
	with $\rho=\frac{r}{\gcd\{r,s\}}$ and $\sigma=\frac{s}{\gcd\{r,s\}}$.
	
\begin{figure}[t]
		\centering
		\includegraphics[width=\linewidth]{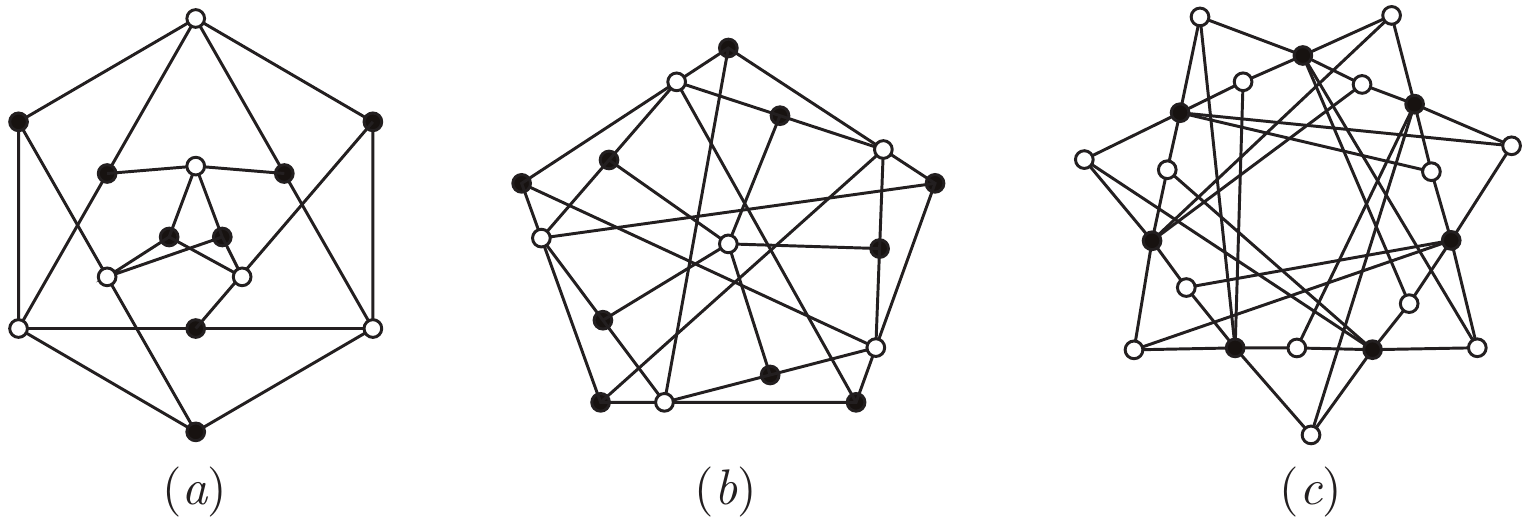}
		\vskip-6.75cm
		\caption{$(a)$ The only [4,3;3]-biMoore graph on $14$ vertices; $(b)$ One of the two [5,3;3]-biMoore graphs on $16$ vertices; $(c)$ The only [6,3;3]-biMoore graph on $21$ vertices.}
		\label{3unics}
	\end{figure}

	Two bipartite biregular graphs with diameter three attaining the Moore bound \eqref{Moore-odd(d=3)}  were given in Yebra, Fiol, and F\`abrega \cite{yff83}. Namely, in Figure \ref{3unics}$(a)$, with $r=4$ and $s=3$, we would have the unattainable values $(N_1',N_2')=(9,10)$, whereas we get $(N_1,N_2)=(6,8)$, giving $M(4,3;3)=14$. In Figure \ref{3unics}$(b)$, with $r=5$ and $s=3$, we have $(N_1',N_2')=(11,13)$, and $(N_1,N_2)=(6,10)$, now corresponding to $M(5,3;3)=16$. In Figure \ref{3unics}$(c)$, with $r=6$ and $s=3$, there is an optimal (unique) $[6,3;3]$-bigraph on 21 vertices, whereas the Moore bound is $M(6,3;3)=8+16=24$. 
	
	In \cite{adfl21}, the authors 
	studied the simple case of biMoore graphs with degrees  $r$ (even), $s=2$, and diameter $d\le 3$. For these values, the Moore bounds in \eqref{Moore-odd(d=3)} turn out to be $M(r,2;3)=2+r$ when $r(>1)$ is odd,  attained by the complete bipartite graph $K_{r,2}$ with $d=2$, and  $M(r,2;3)=3\left(1+\frac{r}{2}\right)$ when $r(>2)$ is even. In the last case, the bound is obtained by a graph with three vertices, $u_1,u_2,u_3$, where each pair $(u_i,u_j)$ is joined by $r$ paths of length 2.
	In the same paper, a numeric construction of bipartite biregular Moore graphs for diameter $d = 3$
	and degrees $r$ and $s=3$ was proposed.
	
	In general, and in terms of designs, to guarantee that the diameter equals $3$, any pair of points must share a block, and any pair of blocks must have a non-empty intersection.
	It is well known that this kind of structure exists when $r=m$ and $r-1$ is a prime power. Namely, the so-called projective plane of order $r-1$. In this case, any pair of blocks (or lines) intersect in exactly one %and only one
	point and, for any pair of points, they share one and only one block (or line).
	For more details about projective planes, you can consult, for instance, Coxeter \cite{Cox93}. This condition is unnecessary in our constructions of block designs, which give graphs of diameter of $3$. Two blocks can be intersected in more than one point, and two points can share one or more blocks. Moreover, we may have the same block appearing more than once. In the following section, we construct graphs with orders that attain or are close to such Moore bounds.
	
	\begin{table}[t]
		\begin{center}
			\begin{tabular}{|c||c|c|c|c|c|c|c|c|c|c|c|}
				\hline
				$r\setminus s$ & 2            & 3   & 4   & 5   & 6   & 7   & 8    & 9   & 10  & 11 & 12 \\
				\hline \hline%\cline{1-2}
				2   &    $\stackrel{(6)}{\textbf{6}}$                                                   \\ \cline{1-3}
				3   & -   & $\stackrel{(14)}{\textbf{14}}$                                               \\ \cline{1-4}
				4   & $\stackrel{(9)}{\textbf{9}}$   & $\stackrel{(14)}{\textbf{14}^{* \diamond}}$       & $\stackrel{(26)}{\textbf{26}}$                                \\ \cline{1-5}
				5   & -   & $\stackrel{(16)}{\textbf{16*}}$            &  27  & $\stackrel{(42)}{\textbf{42}}$   \\ \cline{1-6}
				6   & $\stackrel{(12)}{\textbf{12}}$ &  $\stackrel{(24)}{\textbf{21}^{\diamond}}$ & $\stackrel{(35)}{30}$                       &  44  & $\stackrel{(62)}{\textbf{62}}$                             \\ \cline{1-7}
				7   & -   & $\stackrel{(20)}{\textbf{20*}}$          & 33                       & 48  &  65 &  86                       \\ \cline{1-8}
				8   & $\stackrel{(15)}{\textbf{15}}$ & $\stackrel{(22)}{\textbf{22*}}$         & 42                       & 52  & 70  &  90  & $\stackrel{(114)}{\textbf{114}}$                \\ \cline{1-9}
				9   & - & $\stackrel{(32)}{28^*}$                          & 39                       & 56  & 80  & 96  &  119 & \textbf{146}          \\ \cline{1-10}
				10 &$\stackrel{(18)}{\textbf{18}}$ & $\stackrel{(26)}{\textbf{26*}}$          &  49 & $\stackrel{(69)}{66}$  &  $\stackrel{(88)}{80}$  & 102 & 126 & 152 & $\stackrel{(182)}{\textbf{182}}$  \\ \cline{1-11}
				11 & - & $\stackrel{(28)}{\textbf{28*}}$  & 45 & 64  & 85  & 108 & 133 & 160  & 189  & 222   \\ \cline{1-12}
				12 & $\stackrel{(21)}{\textbf{21}}$ & $\stackrel{(40)}{35^*}$ & 60 & 68 & 99 & 114 & 145 & 175 & 198 & 230 & $\stackrel{(266)}{\textbf{266}}$ \\
				\hline
			\end{tabular}
		\end{center}
		\vskip-.25cm
		\caption{Moore bounds \eqref{Moore-odd(d=3)} for diameter $d=3$ (between parenthesis when a better or optimal construction is known). The known optimal values (most coinciding with the Moore bound) are in boldface. The asterisks correspond to the graphs obtained in \cite{adfl21} 
			(see also Proposition \ref{propo:(r,3)} of that paper), and the diamonds correspond to unique graphs.
			The cases when $r=s=q+1$, with $q$ a prime power, correspond to the point-line incidence graphs of the projective planes.}
		\label{tab:d=3}
		\vskip1cm
	\end{table}

	\subsection{An improved Moore bound for diameter 3}
	\label{sec:Moore-millorada-d=3}
	
	In this subsection, we prove that, for diameter 3 and some values of the degrees $r$ and $s$, the Moore bound in \eqref{Moore-odd(d=3)} can be improved. We begin with the case $r=s^2$.
	
	\begin{lemma}
		Let $G=(V_1\cup V_2,E)$ be an $[r,s;3]$-bigraph, with $r\geq s$, such that $N_2=s(r-1)$. Then, the following conditions hold:
		\begin{itemize}
			\item[$(i)$] 
			Either $r=s$ or $r=s^2$.
			\item[$(ii)$] 
			For any $u \in V_2$, there exists a unique vertex namely $r(u) \in V_2$ (the {\em repeat} of $u$) such that $|N(u)\cap N(r(u))|=2$, where $N(u)$ is the set of vertices that are neighbors of $u$ (at distance $1$).
		\end{itemize}
		\label{lem:newbound3}
	\end{lemma}
	
	\begin{proof}
		$(i)$ From $rN_1=sN_2$ and assuming $N_2=s(r-1)$, we have $N_1=\frac{s^2(r-1)}{r}$. Besides, $\gcd(r,r-1)=1$ and, as a consequence, $r$ must divide $s^2$. Since $r \geq s$, the only possibilities are either $r=s$ or $r=s^2$.\\ 
		$(ii)$ Now, we prove that, for any $u \in V_2$, there exists a unique vertex, namely, $r(u) \in V_2$ (the {\em repeat} of $u$) such that $|N(u)\cap N(r(u))|=2$. There are exactly two different paths from $u$ to $r(u)$ of length $2$. Indeed, since the diameter of $G$ is three, for every $u' \in V_2$, with $u'\neq u$, $\dist(u,u')=2$ must hold. The maximum number of vertices at distance $2$ from $u$, allowed by the vertex degrees, is $s(r-1)$ and since $N_2=s(r-1)$ by hypothesis. Then, the existence of a vertex $r(u)$ sharing two neighbors with $u$ is guaranteed. Moreover, for the remaining vertices $u'$ of $V_2$, others than $r(u)$, there must be a unique path of length $2$ from $u$ to $u'$. 
	\end{proof}
	The case $r=s^2$ is very interesting since the condition $N_2=s(r-1)$ provides the Moore bound in \eqref{Moore-odd(d=3)}. In fact, as a consequence of Lemma \ref{lem:newbound3}, such a bound can be improved. 
	
	\begin{proposition} 
		\label{nonupperbound}
		Given $s\ge 3$, there is no $[s^2,s;3]$-bigraph with order attaining the Moore bound in \eqref{Moore-odd(d=3)}, $M(s^2,s;3)=(s^2-1)+(s^3-s)$.  %, for any $s\geq 3$. 
		Instead, the new, improved Moore bound is
		\begin{equation}
			\label{Moore-odd-improved(r+1,r)}
			M^*(s^2,s;3)= (s^2-2)(s+1),
		\end{equation}
		obtained with
		\begin{equation}
			\label{N1-N2-(r+1,r)a}
			N_1= s^2-2
			\qquad \mbox{and}\qquad N_2=s^3-2s.
		\end{equation}
	\end{proposition}
	
	\begin{proof}
		Assume that $M(s^2,s;3)$ is attained and, hence, $N_1=s^2-1$ and $N_2=s^3-s$. Then, we claim that, for any $v \in V_1$, there exists another vertex $v' \in V_1$ such that $|N(v) \cap N(v')| \geq 3$. Indeed, assume to the contrary that $|N(v)\cap N(w)|\leq 2$ for all $w \in V_1$, with $w\neq v$. Then, the number of edges from $V_1\setminus \{v\}$ to $N(v)$ would be $\leq 2(N_1-1)$. Moreover, $|N(v)|=s^2$ and, since every vertex in $N(v)$ has degree $s$, the number of edges from $N(v)$ to $N^2(v)$$(=V_1\setminus \{v\}$) is $s^2(s-1)$, where $N^2(v)$ is the set of vertices at distance 2 from $v$. As a consequence, $s^2(s-1)\leq 2(N_1-1)$. Hence, 
		$N_1\ge \frac{1}{2}s^2(s-1)+1$, which, for $s\ge 3$, is greater than $s^2-1$, a contradiction. 
		Thus, for any $v \in V_1$, there exists another vertex $v' \in V_1$ such that $|N(v) \cap N(v')| \geq 3$. Let $N(v) \cap N(v')=\{u_1,u_2,\dots,u_l\}\subset V_2$, where $l \geq 3$. Notice that every $u_i$, with $i=2,\dots,l$, is a repeated vertex of $u_1$ since there are two different paths $u_1-v-u_i$ and $u_1-v'-u_i$ joining them, but this is a contradiction with the fact that there is a unique repeated vertex for any vertex in $V_2$, as we proved in Lemma \ref{lem:newbound3}.
		Therefore, our initial assumption is false, and the new possible values of $N_1$ and $N_2$ satisfying $s^2N_1=sN_2$ are those in \eqref{N1-N2-(r+1,r)a}.
		% \textst{The Moore bound given in \eqref{Moore-odd(d=3)} for the $[s^2,s;3]$ case is $M(s^2,s;3)=s^3+s^2-s-1$, where $N_1=s^2-1$ and $N_2=s^3-s$. So, taking $r=s^2$ in Lemma \ref{lem:newbound3}, we have that there is no $[s^2,s;3]$-bigraph $G$ such that $N_2=s(r-1)=s(s^2-1)=s^3-s$ and $N_1<\frac{1}{2}s^2(s-1)+1$. The inequality $s^2-1<\frac{1}{2}s^2(s-1)+1$ holds for any $s \geq 3$.}
	\end{proof}
	
	The above proposition can be generalized to give an improved Moore bound for $[\rho s,s;3]$-bigraphs under some conditions on $\rho$ and $s$. With this aim, we extend the concept of repeat to {\em redundant paths}, defined as expected: The number of redundant paths between a vertex $u$ to a vertex subset $U$  is the total number of shortest paths from $u$ to $U$ minus $|U|$.
	
	\begin{proposition} 
		\label{nonupperbound(b)}
		Given $s\ge 3$ and $\rho$ such that $s-1\le \rho\le s^2-s-3$, there is no $[\rho s,s;3]$-bigraph with order attaining the Moore bound in \eqref{Moore-odd(d=3)}, $M(\rho s,s;3)=(s^2-1)(\rho+1)$. Instead, the new, improved Moore bound is
		\begin{equation}
			\label{Moore-odd-improved(r+1,r)b}
			M^*(\rho s,s;3)= (s^2-2)(\rho+1),
		\end{equation}
		obtained with
		\begin{equation}
			\label{N1-N2-(r+1,r)b}
			N_1= s^2-2
			\qquad \mbox{and}\qquad N_2=\rho(s^2-2).
		\end{equation}
	\end{proposition}

\begin{figure}[t]
		\centering
		\vskip-1cm
		\includegraphics[width=0.8\textwidth]{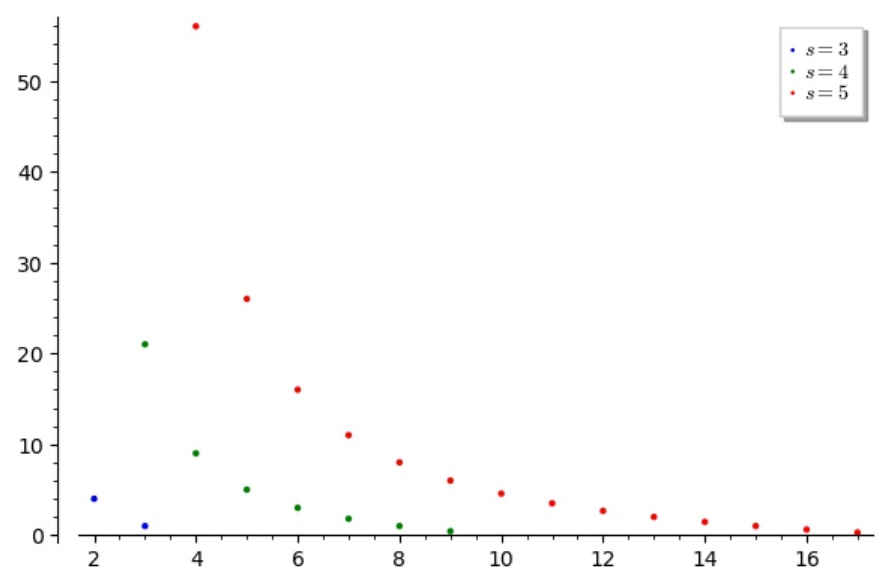}
		\caption{The values of the function $\phi(\rho,s)$ for $s=3$ (blue), $s=4$ (green) and $s=5$ (red).}
		\label{fig:phi}
	\end{figure}

	\begin{proof}
		Suppose that, by contradiction, the bound $M(\rho s,s;3)$ with  $N_1=s^2-1$ and $N_2=\rho(s^2-1)$ is attained. Then, the number of redundant $2$-paths from $u\in V_2$ to $V_2\setminus \{u\}$ is 
		$$
		s(r-1)-(N_2-1)=s(\rho s-1)-\rho(s^2-1)+1=\rho-s+1.
		$$
		%
		%
		% and, similarly,
		% the number of redundant $2$-paths from $v\in V_1$ to $V_1$ is 
		% $$
		% r(s-1)-(N_1-1)=\rho s(s-1)-(s^2-2)=(\rho-1)s^2-\rho+2.
		% $$
		Then, we claim that, for any $v \in V_1$, there exists another vertex $v' \in V_1$ such that $|N(v) \cap N(v')| \geq \rho-s+3$. If not, we would have  that $|N(v)\cap N(w)|\leq \rho-s+2$ for all $w \in V_1$, with $w\neq v$. Hence, the number $e=\rho s(s-1)$ of edges from $V_1\setminus \{v\}$ to $N(v)$ would be at most $(\rho-s+2)(N_1-1)$. 
		% Moreover, $|N(v)|=s^2$ and, since every vertex in $N(v)$ has degree $s$, the number of edges from $N(v)$ to \blue{$N^2(v)(=V_1\setminus \{v\}$)} is $s^2(s-1)$. 
		As a consequence, from $\rho s(s-1)\leq (\rho-s+2)(N_1-1)$, we would have that
		$
		N_1\ge \frac{\rho s(s-1)}{\rho -s +2}.
		$
		However, the function 
		$$
		\phi(\rho,s)= \frac{\rho s(s-1)}{\rho -s +2}-(s^2-1)
		$$ 
		is positive for $s\ge 3$ and $s-2<\rho<s^2-s-2$ (see Figure \ref{fig:phi}), a contradiction. 
		Thus, for any $v \in V_1$, there exists another vertex $v' \in V_1$ such that $|N(v) \cap N(v')| \geq \rho-s+3$. From here, we reason as in the proof of Proposition \ref{nonupperbound} to conclude \eqref{Moore-odd-improved(r+1,r)b}.
		%  Let $N(v) \cap N(v')=\{u_1,u_2,\dots,u_l\}\subset V_2$ where $l \geq 3$. Notice that every $u_i$, $i=2,\dots,l$, is a repeated vertex of $u_1$ since there are two different paths $u_1-v-u_i$ and $u_1-v'-u_i$ joining them, but this is a contradiction with the fact that there is a unique repeated vertex for any vertex in $V_2$, as we proved \blue{in Lemma \ref{lem:newbound3}.
			% Therefore, our initial assumption is false, and the new possible values of $N_1$ and $N_2$ satisfying $s^2N_1=sN_2$ are those in \eqref{N1-N2-(r+1,r)a}.}
	\end{proof}

	See some examples of these improved new values of the Moore bound in Table \ref{tab:d=3(b)}, for $s=3$, 4, and 5.
	
	\begin{table}[t]
		\begin{center}
			\begin{tabular}{|c||c|c|c|c|c|c|c|c|c|c|c|c|}
				\hline
				$s\setminus r$ & 6 & 9 & 12 & 16 & 20  & 24 & 25 & 28 & 30 & 32   &  35  & 36 \\
				\hline \hline%\cline{1-2}
				3   &    {\bf 21}   & {\bf 28} \\ 
				\cline{1-13}
				4   & -   & -   & 56 & 70 & 84 & 98 & - & 112 & - & 126 & - & 140                                            \\ \cline{1-13}
				5   & -   & -  & - & - & 115 & - & 138 & - & 161 &  - & 184 & -\\  
				\hline
			\end{tabular}
		\end{center}
		\vskip-.25cm
		\caption{Improved Moore bounds \eqref{Moore-odd-improved(r+1,r)b} for diameter $d=3$. The known optimal values are in boldface. (Compare with the values in Table \ref{tab:d=3}.) }
		\label{tab:d=3(b)}
		%\vskip1cm
	\end{table}
	
	\section{Bipartite biregular graphs derived from difference sets}
	\label{sec:d=3}
	
	Given a group $(\Gamma,\cdot)$ of order $n$, an $(n,s,\lambda)$-difference set in $\Gamma$ is an $s$-subset $S$ of $\Gamma$ such that each non-identity element $g \in \Gamma$ occurs exactly $\lambda$ times in the multiset $\{t_i\cdot t_j^{-1} \, | \, t_i,t_j \in S\}$. The study of the difference sets of groups dates back to the thirties of the past century (see Singer \cite{js38}). His study is focused mainly on the difference sets of Abelian groups, where the additive notation $\{t_i-t_j \, | \, t_i,t_j \in S\}$ is commonly used. In particular, a difference set with $\lambda=1$ is called  {\em perfect}. For instance, in the case of $\Gamma=\Z_{13}$, we can take  $S=\{0,1,3,9\}$. As shown in the following table, here in additive notation, $d_{ij}=i-j\pmod{13}$ for $i,j\in S$, the {\em difference matrix $\D_S=(d_{ij})$} gives each non-zero residue modulo 13 exactly once. That is, $S$ is a $(13,4,1)$-difference set.
	$$
	\begin{tabular}{|c||c|c|c|c|}
		\hline
		$i\backslash j$ & 0 & 1 & 3 & 9\\
		\hline
		\hline
		0 & 0 & 12 & 10 & 4\\
		\hline
		1 & 1 & 0 & 11 & 5\\
		\hline
		3 & 3 & 2 & 0 & 7\\
		\hline
		9 & 9 & 8 & 6 & 0\\
		\hline
	\end{tabular}\ , \qquad \D_S=
	\left(\begin{array}{cccc}
		0 & 12 & 10 & 4 \\
		1 & 0 & 11 & 5 \\
		3 & 2 & 0 & 7 \\
		9 & 8 & 6 & 0
	\end{array}\right).
	$$
	Let $S$ be a subset of cardinality $s$ of a given group $\Gamma$ of order $n$. We say that $S$ is a {\em covering difference set} of $\Gamma$ if the multiset of differences $\{t_i\cdot t_j^{-1}\, | \, t_i,t_j\in S\}$ contains all the elements of $\Gamma$ (perhaps with some repetitions). That is, the difference matrix $\D_S$ covers all the elements of $\Gamma$. We use these sets to construct bipartite biregular graphs with large order as follows.
	
	\begin{definition}
		Given $\Gamma$ and $S$ as above, and an integer $m \geq 1$, let $G_m(S)=(V_0 \cup V_1, E)$ be the bipartite graph with independent sets $V_0=\{(0,u): u\in \Gamma\}$, $V_1=\{ (l,v): v\in \Gamma\ \textrm{ and }\ 1 \leq l \leq m\}$, and adjacencies $(l,v) \sim (0,v\cdot \sigma)$ for all $\sigma\in S$ and $l \in \{1,\ldots,m\}$. 
	\end{definition}
	
	The next result gives the main properties of these graphs.
	
	\begin{proposition}
		\label{propo:(r,3)}
		Let $\Gamma$ be an Abelian group with $n$ elements and a covering $s$-set $S$. Then,  the bipartite graph $G_m(S)$ has  $(m+1)n$ vertices, where every vertex in $V_0$ has degree $s$ and every vertex in $V_1$ has degree $ms$. Moreover, the diameter of $G_m(S)$ is $3$.
	\end{proposition}
	
	\begin{proof}
		The order of $G_m(S)$ follows directly from the definition: $V_0$ contain the $n$ elements of $\Gamma$, meanwhile $V_1$ contains $m$ different copies of $\Gamma$. Let $S=\{t_1,\dots,t_s\}$. Notice that each vertex $(j,v) \in V_1$ is adjacent to $(0,v\cdot t_1),(0,v\cdot t_2),\dots,(0,v\cdot t_s)$ and all these vertices are different since $\Gamma$ is a group. Besides, each vertex of $V_2$, say $(0,v)$, is adjacent to $(l,v\cdot t_i^{-1})$ for all $1 \leq l \leq m$ and $t_i \in S$. This gives the desired degrees. To see that the diameter is $3$, it suffices to prove that every pair of different vertices $(0,u),(0,u')\in V_1$, and every pair of different vertices $(l,v),(l',v')\in V_1$, are at distance two.
		
		In the second case, we have the following path of length $2$:
		$$
		(l,v)\  \sim \  (0, v\cdot t_i)\  \sim \  (l',v\cdot t_i \cdot t_j^{-1}) \quad \mbox{for all $t_i,t_j \in S$ and $l,l' \in \{1,\ldots,m\}$}\nonumber.
		$$
		
		Notice that $t_i \cdot t_j^{-1}$ is precisely the element $d_{ij}$ in the difference matrix $\D_S$ and, hence, it provides every element of $\Gamma$ under the hypothesis. A path of length $2$ between vertices in $V_0$ is given by:
		\begin{equation}\label{eq:graph}
			(0,v)\  \sim \  (l, v\cdot t_i^{-1})\  \sim \  (0,v\cdot t_i^{-1} \cdot t_j) \quad \mbox{for all $t_i,t_j \in S$ and $l \in \{1,\cdots,m\}$}.   
		\end{equation}
		Since $\Gamma$ is an Abelian group, then  $t_i^{-1} \cdot t_j =t_j\cdot t_i^{-1}$, and so every element in $\Gamma$ can be represented as $v\cdot t_i^{-1} \cdot t_j$ by the hypothesis.
	\end{proof}
	
	\begin{figure}[t]
		\centering
		\includegraphics[width=6cm]{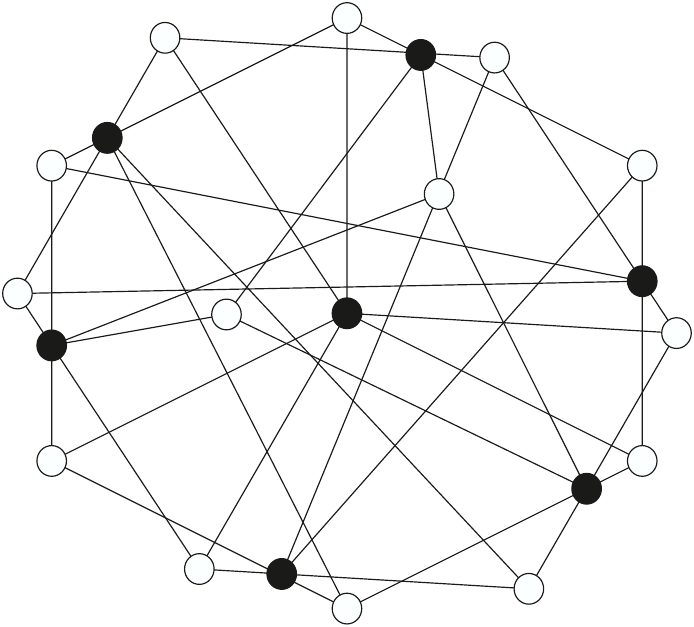}
		%\vskip-6.75cm
		\caption{The graph $G_2(\{0,1,3\})$ over $\mathbb{Z}_7$. It has $21$ vertices, degrees $(6,3)$ and diameter 3.}
		\label{G(7,14)}
	\end{figure}
	
	For example, in Figure \ref{G(7,14)}, there is the graph $G_2(S)$ obtained from the perfect set $S=\{0,1,3\}$ in $\Gamma=\mathbb{Z}_7$. It has 21 vertices, degrees $(6,3)$, and diameter $d=3$. Although the corresponding Moore bound is 24, this graph is known to be optimal. In fact, $S$ is a perfect difference set from where the best results are obtained.
	%\blue{Recall that a {\em difference set} for some group. 
		Another example is the graph $G_2(S)$ of Figure \ref{G_2(0,1,3,9)}, where $S=\{0,1,3,9\}$ is the perfect difference set on $\Z_{13}$ shown at the beginning of this section. Then, it has 39 vertices, degrees $(8,4)$, and diameter 3.
		
		A simple counting argument says that any $(n,s,\lambda)$-difference set $S$ must satisfy $\lambda(n-1)=s(s-1)$. Given an $s$-subset $S$, the equality $n=s^2-s+1$ must hold to look for a perfect difference set in a group $\Gamma$ of order $n$. 
		
		Then, it is known that, for every prime power $q$, there exists a perfect  $(n,s,1)$-difference set $S$ with $s=|S|=1+q$ elements, and $n=q^2+q+1=s^2-s+1$ (see Singer \cite{js38}, and Halberstam and Laxton \cite{hl63}). 
		In fact, Singer proved that if
		$S=\{s_i:0\le i\le q\}$ is a perfect difference set modulo $n=q^2+q+1$, and if lines $L_i$ and points $P_j$ (for $i,j\in\Z_n$) of $S$ are defined in such a way that each line $L_i$ contains the set of points $P_j=i+s_j$ for  $j=0,1,\ldots, q$, then $S$ is a projective plane.
		As a consequence, the graphs $G_1(S)$ of Proposition \ref{propo:(r,3)} with $r=s=q+1$ turn out to be the point-line incidence graphs of the projective planes $PG(2,q)$.
		
		A handy way of obtaining a difference set $S$ with
		$s=|S|=q+1$ elements, relative to the group $\Z_{n}$ with $n=q^2+q+1$, is as follows. 
		Let $\phi(x)$ be a polynomial of degree $3$, with coefficients in $GF(q)\cong \Z_q$, irreducible (so that $\phi(\alpha)\neq 0\ (\textrm{mod }q)$ for every $\alpha\in GF(q)$) and primitive (so that $r=q^n-1$ is the smallest integer such that $\phi(x)$ divides $x^r-1$).
		Then, in the extension field $F=GF(q^3)$ of polynomials $\textrm{mod }\phi(x)$, the monic polynomial $x$ is primitive. Hence, the elements of the multiplicative (cyclic) group $F^*$ of $GF(q^3)$ can be represented as the powers $x^0=1, x^1=x, x^2,\ldots, x^{r-1} (\textrm{mod }\phi(x))$. Now, suppose that, in $F$, the elements $1+x,1+2x,\ldots, 1+(q-1)x$ correspond to $x^{s_1},x^{s_2},\ldots, x^{s_{q-1}}$, respectively. Then, the set $S=\{0,1,s_1,s_2,\ldots,s_{q-1}\}$ is a difference set with $q+1$ elements, relative to $\Z_{q^2+q+1}$. (Of course, all the numbers $s_i$, for $i=1,2,\ldots, q-1$, can be reduced modulo $m$ so that $2\le s_i\le n-1$).
		For example, if $q=3$, a primitive polynomial to derive the extension field $F=GF(3^3)$ is $\phi(x)=x^3+2x^2+x+1$. Then, the elements of $F$ are:
		\begin{center}
			\begin{tabular}{llll}
				$x^0=1$, & $x^1=x$ & $x^2=x^2$ & $x^3=x^2+2x+2$ \\
				$x^4=x+2$, & $x^5=x^2+2x$, & $x^6=2x+2$, & $x^7=2x^2+2x$, \\ $x^8=x^2+x+1$, & $x^9=2x^2+2$ & $x^{10}=2x^2+1$, & $x^{11}=2x^2+2x+1$\\
				$x^{12}=x^2+2x+1$, & $x^{13}=2$, & $x^{14}=2x$, & $x^{15}=2x^2$,\\
				$x^{16}=2x^2+x+1$, & \boldmath{$x^{17}=2x+1$}, & $x^{18}=2x^2+x$, & \boldmath{$x^{19}=x+1$},\\
				$x^{20}=x^2+x$, &
				$x^{21}=2x^2+2x+2$, &  $x^{22}=x^2+1$, & $x^{23}=x^2+2$, \\ $x^{24}=x^2+x+2$, & $x^{25}=2x^2+x+2$. & &
			\end{tabular}
		\end{center}
		
		Therefore, a difference set for $\Z_n$, with $n=q^2+q+1=13$, is
		$S=\{0,1,17,19\ (\textrm{mod }13)\}=\{0,1,4,6\}$ (see the table below).
		$$
		\begin{tabular}{|c||c|c|c|c|}
			\hline
			$i\backslash j$ & 0 & 1 & 4 & 6\\
			\hline
			\hline
			0 & 0 & 12 & 9 & 7\\
			\hline
			1 & 1 & 0 &  10 & 8\\
			\hline
			4 & 4 & 3 & 0 & 11\\
			\hline
			6 & 6 & 5 & 2 & 0\\
			\hline
		\end{tabular}
		$$
		
		We get the following consequence from Proposition \ref{propo:(r,3)}.
		\begin{corollary}
			\label{coro:perfect}
			Let $s=1+q$ with $q$ a prime power. Let $n=s^2-s+1$. Then, for every $m\ge 1$, there is a bipartite graph $G_m(S)$ on $(m+1)n$ vertices and degrees $r=ms$ and $s$ with diameter $d=3$.
		\end{corollary} 
		
		Thus, these graphs $G_m(S)$, with pair of degrees $(ms,s)$, have order 
		\begin{equation}
			|V_0|+|V_1|=n(m+1)=(s^2-s+1)(m+1),
			\label{order}
		\end{equation}
		whereas, the Moore bound \eqref{Moore-odd(d=3)} for degrees $(r,s)=(ms,s)$ and diameter $3$ is 
		\begin{equation}
			M(ms,s;3)=(m+1)\left\lfloor s^2-\frac{s-1}{m}\right\rfloor = (s^2-1)(m+1)
			\label{moore-d=3}
		\end{equation}
		when $m\ge s+1$, and
		the improved Moore bound \eqref{Moore-odd-improved(r+1,r)b} is
		\begin{equation}
			M^*(ms,s;3)= (s^2-2)(m+1)
			\label{fita-Moore-millorada}
		\end{equation}
		when $s-1\le m\le s^2-s-3$.
		% \red{Oigan pero justamente nosotros acabamos de mejorar estas cotas con la Proposici\'on 1.3 para valores $(ms,s;3)$
			% $$M^*(ms,s;3)= (s^2-2)(m+1)$$
			% Por ejemplo, en la tabla que aparece a continuación, acabamos de probar, que para $s=3$, tenemos $(3m,3;3)$-biMoore grafos, es decir las que construimos alcanzan la "nueva" cota superior que ahora es $7(m+1)$ y para $s=4$, es decir para $(4m,4;3)$-bigrafos, el grafo es, por construcci\'on de orden $13(m+1)$ pero la cota superior es ahora $14(m+1)$ Y bueno insisto que a mi me gustaría probar, ya sea con la existencia de los planos proyectivos o los conjuntos diferencia (que son lo mismo) que podríamos "reducir más" esa cota.}
		%(with equality when $m\geq s-1$). 
		So, whenever the graphs are in the setting of Corollary \ref{coro:perfect},  their order asymptotically approaches the Moore bound \eqref{moore-d=3} of diameter $3$ with the ratio
		$$
		\gamma= \frac{s^2-s+1}{s^2-1}.
		$$
		In particular, when $s-1\le m\le s^2-s-3$, the values in \eqref{order} and \eqref{fita-Moore-millorada} coincide, and we get $(3m,3;3)$-biMoore graphs.
		% for degrees $(ks,s)$, where $k \geq s-1$. The results are even better when $k<s$ since the Moore bound is smaller.
		In Table \ref{tab:singer}, we show perfect difference sets for $s\le 12$, together with the parameters
		of the corresponding graphs $G_m(S)$ and the improved Moore bound \eqref{Moore-odd-improved(r+1,r)b} provided that $s-1\le m\le s^2-s-3$ (otherwise, we must use \eqref{Moore-odd(d=3)}). For instance, the family of graphs $G_m(S)$, constructed with the $(13,4,1)$-difference set 
		$S=\{0,1,3,9\}$, produces a bipartite graph with degrees $r=4m$ and $s=4$ having $13(m+1)$ vertices and diameter $d=3$.
		
		\begin{figure}[t]
			\centering
			\includegraphics[width=8cm]{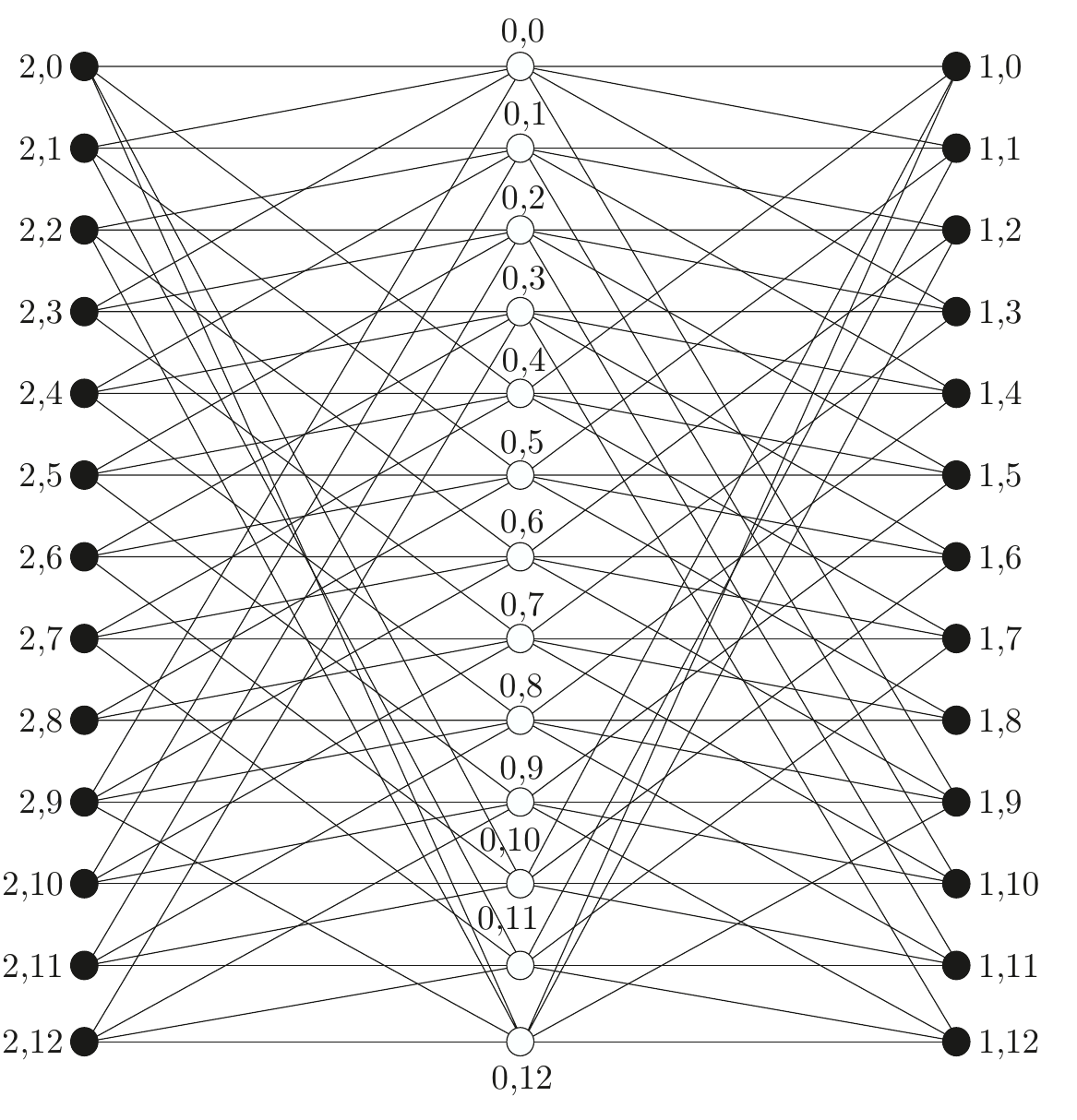}
			\caption{The bipartite graph $G_2(\{0,1,3,9\})$ over $\mathbb{Z}_{13}$. It has $39$ vertices, degrees $(8,4)$, and diameter 3.}
			\label{G_2(0,1,3,9)}
		\end{figure}
		
		The corresponding Moore bound for this case $(4m,4;3)$ is $42$ for $m=2$, and $15(m+1)$ for $m\geq 3$. As a consequence, the family of graphs $G_m(S)$ asymptotically approaches the Moore bound with a constant ratio of $\frac{13}{15} \approx 0.867$. In particular, for $m=2$, the corresponding graph with degrees $(r,s)=(8,4)$ has 39 vertices (see Figure \ref{G(7,14)}), whereas the Moore bound is 42 (see Table \ref{tab:d=3}). Another good example is taking the set $S=\{0,1,4,14,16\}$, which is a perfect difference set for $\Z_{21}$. The Moore bound for the set $(5m,5;3)$ is $24(m+1)$ (for $m \geq 5$), and it provides $21(m+1)$ vertices. This graph asymptotically approaches the Moore bound in $\frac{7}{8}=0.875$. For $m=2$, we obtain a graph with degrees $(r,s)=(10,5)$ having 63 vertices, whereas the maximum order is 66. 
		
		\begin{table}[t]
			\begin{center}
				\begin{tabular}{|c|c|c|c|c|}
					\hline
					$s$ & $\Gamma$ & $S$ & order of    & Moore \\
					&          &     & $G_m(S)$ & Bound \eqref{Moore-odd-improved(r+1,r)b} \\
					\hline
					\hline
					3 & $\mathbb{Z}_7$ & $\{0,1,3\}$ & $7(m+1)$ & $7(m+1)$\\
					\hline
					4 & $\mathbb{Z}_{13}$ & $\{0,1,3,9\}$ & $13(m+1)$ & $14(m+1)$\\
					\hline
					5 & $\mathbb{Z}_{21}$ & $\{ 0, 1, 4,14,16\}$ & $21(m+1)$ & $23(m+1)$\\
					\hline
					6 & $\mathbb{Z}_{31}$ &$\{ 0, 1, 6,18,22,29\}$ & $31(m+1)$ & $34(m+1)$ \\
					\hline
					8 & $\mathbb{Z}_{57}$ & $\{ 0, 1, 5, 7,17,35,38,49\}$ & $57(m+1)$ & $62(m+1)$\\
					\hline
					9 & $\mathbb{Z}_{73}$ & $\{ 0, 1,17,39,41,44,48,54,62\}$ &$73(m+1)$ & $79(m+1)$ \\
					\hline
					10 & $\mathbb{Z}_{91}$ & $\{0, 1, 3, 9,27,49,56,61,77,81\}$ & $91(m+1)$ & $98(m+1)$\\
					\hline
					12 & $\mathbb{Z}_{133}$ & $\{0, 1, 3, 12, 20, 34, 38, 81, 88, 94, 104, 109\}$ & $133(m+1)$ & $142(m+1)$\\
					\hline
				\end{tabular}
			\end{center}
			\caption{Perfect difference sets (Singer \cite{js38}) for some values of $s=|S|$ and the parameters of the corresponding bipartite graphs of Corollary \ref{coro:perfect}, with degrees $(ms,s)$, diameter $3$, and the Moore bounds \eqref{Moore-odd-improved(r+1,r)b} (for $s-1\le m\le s^2-s-3$).}
			\label{tab:singer}
		\end{table}
		
		It is conjectured that perfect difference sets only exist whenever $s-1$ is a prime power. The conjecture has been proved true for $s-1 \leq 1600$ (see Evans and Mann \cite{evansmann51}), but an answer to the general question remains an open problem. 
		In all those cases in which perfect difference sets do not exist, it would be useful to get other difference sets that are `close' to being perfect, in the sense that every element of the group appears as a difference of two elements of $S$, but maybe this representation is not unique. There is an approximation to difference sets in the literature that fits in this way. An $(n, s, \lambda, t)$-{\em almost difference set} (ADS) is a subset $S$ of size $s$ taken from an order-$n$ group such that the multiset of differences contains $t$ of the nonidentity group elements $\lambda$ times and all other nonidentity group elements $\lambda + 1$ times. The problems of existence and construction for ADS's in arbitrary groups are even less well understood than for general difference sets. The significantly weaker difference property required, along with the relatively recent development of the ADS concept, implies that the theory of ADS's is still a work in progress. 
		
		Inspired by the above results, we end this section by posing the following (ambitious) conjecture.
		\begin{conjecture}
			For $s-1$ a prime power, the graphs $G_m(S)$ of Corollary \ref{coro:perfect} are $(ms,s;3)$-biMoore graphs
		\end{conjecture}
		As mentioned, the conjecture holds when $m=1$, where $G_1(S)$ is the point-line incidence graph of $PG(2,s-1)$.

		\section{Covering difference sets for \textit{s}=7}
		\label{sec:d=4}
		There are many values of $s$ in which perfect difference sets $S$ do not exist. The first one is $s=7$, and it is interesting on its own. 
		The order $n$ for an Abelian group $\Gamma$ having a covering difference set of order $s=7$ is $n \leq 43$. Such a covering difference set for $n=43$ ($\Gamma=\mathbb{Z}_{43})$ would produce a perfect difference set, but we already know that it does not exist (see Evans and Mann  \cite{evansmann51}). So, we should look for covering sets of size $7$ in Abelian groups of order $n\leq 42$. We have done an exhaustive computer search for finding a cover difference set in $\Gamma=\mathbb{Z}_{42}$, but we found none. To this end, notice that there are $\binom{42}{7} \approx 27$ million candidates to test, but two properties reduce the computational cost of the exhaustive search. The first property is that $S=\{t_1,t_2,\dots,t_s\}$ is a covering difference set if and only if $S-t_1=\{0,t_2-t_1,\dots,t_s-t_1\}$ is a covering difference set. In our case, the proof is trivial from the definition and reduces the number of testing sets to $\binom{41}{6} \approx 4.5$ million. The second property is the following result. 
		
		\begin{proposition}\label{prop:rep}
			Let $\Gamma$ be a group of order $n$ with a covering set $S$ of order $s$, such as $n=s^2-s$. Then, every $g \in \Gamma\setminus \{0\}$ appears once in $D_S$ except for just one element of order $2$ that appears twice.
		\end{proposition}
		\begin{proof}
			Notice that $D_S$ has $s^2-s$ elements out of the main diagonal and, hence,  each $g \in \Gamma\setminus \{0\}$ must appear once in $D_S$ except one element $g'$ appearing twice (to fit the $s^2-s-1$ different elements of $\Gamma\setminus \{0\}$ out of the diagonal of $D_S$). Moreover, since $g'=t_i-t_j=t'_i-t'_j$ for different $t_i,t_j,t'_i,t'_j \in S$, then $g'$ must have order $2$ since otherwise its inverse (opposite) $-g'=t_j-t_i=t'_j-t'_i$ would appear also twice in $D_S$.   
		\end{proof}
		
		This property allows us to exclude many sets very quickly, namely, all the elements different from $21$ that appear more than once during the computational building of $\D_S$. Then, the set is automatically excluded. That is, we do not need to complete the calculation of $\D_S$ to exclude many candidate sets.
		Of course, Proposition \ref{prop:rep} can be extended to other values of $n$ that are `close' to the maximum $s^2-s+1$ by allowing a maximum number of repetitions for an element to be in $\D_S$ when checking if $S$ can be a covering set of $\Gamma$. This has been very useful in performing an exhaustive search of covering sets for every Abelian group of orders $42$, $41$, and $40$. In every case, we found none. In the first case, $n \leq 42$, a covering set exists when $n=39$.
		
		\begin{proposition}
			The maximum order $n$ for an Abelian group $\Gamma$ having a covering set $S$ of order $s=7$ is $n=39$. The set $S=\{0, 1, 2, 4, 13, 18, 33\}$ in $\mathbb{Z}_{39}$ provides a $(39,7,1,34)$ almost difference set (ADS).
		\end{proposition}
		
		The set $S=\{0, 1, 2, 4, 13, 18, 33\}$ in $\mathbb{Z}_{39}$ provides the following difference matrix $\D_S=(d_{ij})$, where $d_{ij}=t_i-t_j\ (\textrm{mod }39)$, for $t_i,t_j \in S$:
		
		$$
		\D_S=
		\left(\begin{array}{rrrrrrr}
			0 & 38 & 37 & 35 & 26 & 21 & 6 \\
			1 & 0 & 38 & 36 & 27 & 22 & 7 \\
			2 & 1 & 0 & 37 & 28 & 23 & 8 \\
			4 & 3 & 2 & 0 & 30 & 25 & 10 \\
			13 & 12 & 11 & 9 & 0 & 34 & 19 \\
			18 & 17 & 16 & 14 & 5 & 0 & 24 \\
			33 & 32 & 31 & 29 & 20 & 15 & 0
		\end{array}\right),
		$$
		where every non-zero element appears just once except for $1,2,37,38$, which appears twice. Hence, this set $S$ is a $(39,7,1,34)$ ADS. Then, the corresponding graphs $G_m(S)$, with degrees $(7m,7)$ and diameter $3$, have $39(m+1)$ vertices (while the Moore bound is $48(m+1)$).

		Nevertheless, a bit more can be done using non-Abelian groups, as shown in the next subsection.
		
		\subsection{Bipartite biregular graphs derived from difference sets of non-Abelian groups}
		\label{sec:difference-sets-non-Abelian}
		
		From the first contributions to the study of difference sets in $\mathbb{Z}_m$ almost one century ago, many published papers have been devoted to difference sets in Abelian groups. Curiously enough, the situation for non-Abelian groups is quite different, where only a few contributions can be found, many of them from the 90's decade of the past century.
		
		Although all the previous definitions related to difference sets are valid for non-Abelian groups, Proposition \ref{propo:(r,3)} cannot be applied to non-Abelian groups. The main reason is the adjacencies in \eqref{eq:graph}, where the Abelian property is necessary to prove that there is a path of length two between two distinct vertices in $V_0$. One way to guarantee the existence of such a path is to impose the extra condition that the set of inverses of $S$, $\overline{S}=\{t_1^{-1},t_2^{-1},\dots,t_s^{-1}\}$, is also a covering set. This is indeed true for Abelian groups, where if $g\in \Gamma$ can be represented as a difference of two elements of $S$, that is, if $g=t_i\cdot t_j^{-1}$, then  $g=t_j^{-1}\cdot(t_i^{-1})^{-1}$ due to the Abelian property. Hence,  $\overline{S}$ is also a covering set.
		
		\begin{corollary}
			\label{cor:(r,3)}
			Let $\Gamma$ be a non-Abelian group. Then, $G_m(S)$ is a bipartite graph with $(m+1)n$ vertices, where every vertex in $V_0$ has degree $s$ and every vertex in $V_1$ has degree $ms$. Moreover, the diameter of $G_m(S)$ is $3$ if both $S$ and $\overline{S}$ are covering sets of $\Gamma$.
		\end{corollary}
		\begin{proof}
			Follow the arguments of the proof of Proposition \ref{propo:(r,3)} and note that, due to the adjacencies in \eqref{eq:graph}, we can write $(0,v\cdot t_i^{-1} \cdot t_j)=(0,v\cdot t_i^{-1} \cdot( t_j^{-1})^{-1}$). Since $\overline{S}$ is a covering set, a path of length $2$ exists between every pair of different vertices of $V_0$.   
		\end{proof}
		
		Using an exhaustive computer search, we have proved that there is no covering set $S$ of size $s=7$ for the five non-Abelian groups of order $42$. The largest order for a non-Abelian group with a covering set is $40$, and although there are $11$ groups of that order, only one has covering sets of $7$. This group is a semidirect product of $\mathbb{Z}_5$ and $\mathbb{Z}_8$ via a square map. It has the following representation as a quotient of a free group:
		$$
		\Gamma_1=\langle a,b \mid a^5 = b^8 = 1, bab^{-1} = a^2 \rangle.
		$$
		\begin{proposition}
			The maximum order $n$ for a non-Abelian group having a covering set $S$ of order $s=7$ is $n=40$. The set $S=\{1,b,b^4,ba,ba^{-1}b^2,ab^{-1},bab^2\}$ in $\Gamma_1$ is a $(40,7,1,36)$ almost difference set (ADS). 
		\end{proposition}
		The difference matrix 
		\begin{align*}
			& \D_s=\\
			& \left(\begin{array}{ccccccc}
				1& b^{-1}& b^{-4}& a^{-1}  b^{-1}& b^{-2}  a  b^{-1}& b  a^{-1}& b^{-2}  a^{-1}  b^{-1}\\
				b& 1& b^{-3}& b  a^{-1}  b^{-1}& b^{-1}  a  b^{-1}& b^{2}  a^{-1}& b^{-1}  a^{-1}  b^{-1}\\ 
				b^{4}& b^{3}& 1& b^{4}  a^{-1}  b^{-1}& b^{2}  a  b^{-1}& b^{5}  a^{-1}& b^{2}  a^{-1}  b^{-1}\\
				b  a& b  a  b^{-1}& b  a  b^{-4}& 1& b  a  b^{-2}  a  b^{-1}& b  a  b  a^{-1}& b  a  b^{-2}  a^{-1}  b^{-1}\\ 
				b  a^{-1}  b^{2}& b  a^{-1}  b& b  a^{-1}  b^{-2}& b  a^{-1}  b^{2}  a^{-1}  b^{-1}& 1& b  a^{-1}  b^{3}  a^{-1}& b  a^{-2}  b^{-1}\\ 
				a  b^{-1}& a  b^{-2}& a  b^{-5}& a  b^{-1}  a^{-1}  b^{-1}& a  b^{-3}  a  b^{-1}& 1& a  b^{-3}  a^{-1}  b^{-1}\\ 
				b  a  b^{2}& b  a  b& b  a  b^{-2}& b  a  b^{2}  a^{-1}  b^{-1}& b  a^{2}  b^{-1}& b  a  b^{3}  a^{-1}& 1\\
			\end{array}\right)
		\end{align*}
		provides every element of $\Gamma_1\setminus\{1\}$ just once, except for $b^{-1}ab^{-1},b^4$ and $ba^{-1}b$ that appear twice. Note that $b^4$ has order $2$, so $b^4=b^{-4}$, and $b^{-1}ab^{-1}$ is the inverse of $ba^{-1}b$. Moreover, it can be check that $\overline{S}=\{1, b^{-1}, b^4, a^{-1
		}b^{-1}, b^{-2}ab^{-1}, ba^{-1}, b^{-2}a^{-1}b^{-1}\}$ is also a covering set. According to Corollary $\ref{cor:(r,3)}$, the corresponding graphs $G_m(S)$, with degrees $(7m,7)$ and diameter $3$, have $40(m+1)$ vertices, improving the bound obtained with the Abelian group $\mathbb{Z}_{39}$.

\section*{Conflict of Interest Statement}	

The authors have no conflict of interest.


\begin{thebibliography}{10}
				
				\bibitem{adfl21}
				G. Araujo-Pardo, C. Dalf\'o, M. A. Fiol, and N. L\'opez,
				Bipartite biregular Moore graphs, 
				{\em Discrete Math.} {\bf 344} (2021) 112582.
				
				% \bibitem{AraJajRam}
				% G. Araujo-Pardo, G. R. Jajcay, and A. Ramos Rivera,
				% On a relation between bipartite biregular cages, block designs and generalized polygons,
				% submitted (2020).
				
				% \bibitem{BamBisRoy}
				% J. Bamberg, A. Bishnoi, and G. F. Royle,
				% On regular induced subgraphs of generalized polygons,
				% {\em J. Comb. Theory  Ser. A} {\bf 158} (2018) 254--275.
				
				% \bibitem{BeBoPaPe83}
				% J. Bermond, J. Bond, M. Paoli, and C. Peyrat,
				% Graphs and interconnection networks: Diameter and Vulnerability,
				% in E. Lloyd (Ed.), \emph{Surveys in Combinatorics: Invited Papers for the Ninth British Combinatorial Conference 1983}
				% (London Mathematical Society Lecture Note Series, pp. 1--30), Cambridge University Press, Cambridge (1983).
				
				%\bibitem{BigIto80}
				%N. L. Biggs and T. Ito.
				%Graphs with even girth and small excess.
				%{\em Math. Proc. Camb. Philos. Soc.} {\bf  88} (1980) 1--10.
				
				\bibitem{b66}
				W. G. Brown, On graphs that do not contain a Thomsen graph, {\em Canad. Math. Bull.} {\bf 9}
				(1966), no. 5,  281--285.
				
				
				%\bibitem{ColDin}
				%C. J. Colbourn and J. Dinitz (eds.),
				%{\em The CRC Handbook of Combinatorial Designs. 2nd ed.},
				%Discrete Mathematics and its Applications, Boca Raton, FL: Chapman \& Hall/CRC (2007) 984 pp.
				
				%\bibitem{ConExo&Jaj}
				%M.\ Conder, G.\ Exoo and R.\ Jajcay. On the limitations of the use of solvable groups in Cayley graph cage constructions.
				%{\em Europ. J. of Combinatorics } {\bf 31} (2010) 1819--1828.
				
				\bibitem{Cox93}
				H. S. M. Coxeter, 
				{\em The Real Projective Plane}, 
				3rd edition, Springer-Verlag, New York, 1993.
				
				% \bibitem{cb92}
				% H. S. M. Coxeter and G. Beck, 
				% {\em The Real Projective Plane}, 
				% Springer-Verlag, Berlin, Heidelberg, 1992.
				
				% \bibitem{Cv75}
				% D. Cvetkovi\'{c},
				% Spectra of graphs formed by some unary operations,
				% \emph{Publ. I. Math.} \textbf{19(33)} (1975) 37--41.
				
				%\bibitem{d94}
				%C. Delorme,
				%Distance biregular bipartite graphs,
				%{\em European J. Combin.} {\bf 15} (1994), no. 3, 223--238.
				
				%\bibitem{ErdSachs63} P. Erd\H{o}s and H. Sachs,
				%Regul\" are Graphen gegebener Taillenweite mit minimaler Knotenzahl,
				%{\em  Wiss. Z. Uni. Halle (Math. Nat.) } {\bf 12} (1963) 251--257.
				
				\bibitem{efs21}
				G. Erskine, P. Fratri\v{c}, and J. \v{S}ir\'a\v{n},
				Graphs derived from perfect difference sets,
				{\em Australas. J. Combin. \textbf{80} (2021), no. 1, 48--56.}
				
				\bibitem{evansmann51}
				T. A. Evans and H. B. Mann, On simple difference sets, \textit{Indian J. Stat.} \textbf{11(3/4)} (1951) 357--364.
				% \bibitem{ExooJaj08}
				% G. Exoo and R. Jajcay,
				% Dynamic cage survey,
				% {\em Electron. J. Combin.}, Dynamic Survey {\bf 16} (2008).
				
				% \bibitem{fh64}
				% W. Feit and G. Higman,
				% The nonexistence of certain generalized polygons,
				% {\em J. Algebra} {\bf 1} (1964) 114--131.
				
				\bibitem{FilRamRivJaj19}
				S. Filipovski, A. Ramos Rivera, and R. Jajcay,
				On biregular bipartite graphs of small excess,
				{\em Discrete Math.} {\bf 342} (2019) 2066--2076.
				
				%\bibitem{A}  G.\ Exoo, R.\ Jajcay, M.\ Ma\vv  caj and J.\ \vv  Sir\' a\vv  n. \newblock
				%On the defect of vertex-transitive graphs of given degree and diameter. \newblock
				%{\em J. Comb. Theory, Series B} {\bf 134} (2019) 322-340.
				
				%\bibitem{ExoJaj&Sir}
				%G.\ Exoo, R.\ Jajcay and J.\ \vv  Sir\' a\vv  n. Cayley cages.
				%{\em J. Algebr. Comb. Volume 38, Issue } {\bf 1} (2013) 209--224.
				
				%{\color{red} \bibitem{GR01}
					%C. Godsil, G. Royle
					%Algebraic Graph Theory.
					%2001 Springer-Verlag, New York, Inc.}
				
				%\bibitem{Fil17}
				%S.\ Filipovski,
				%On bipartite cages of excess $4$,
				%{\em Electron. J. Combin.} {\bf 24(1)} (2017) \#P1.40.
				
				%\bibitem{Fil&Jaj}
				%S.\ Filipovski and R.\ Jajcay. \newblock
				%On the excess of vertex-transitive graphs of given degree and girth. \newblock
				%{\em Discrete Math.} {\bf 341}, No. 3 (2018) 772--780.
				
				%\bibitem{Fri71}
				%H.\ D.\ Friedman. \newblock
				%On the impossibility of certain Moore graphs. \newblock
				%{\em J. Comb. Theory } {\bf 10} (1971) 245--252.
				
				%\bibitem{FuLaSeUsWol95}
				%Z.\ Furedi, F.\ Lazebnik, A.\ Seress, V.\ A.\ Ustimenko and A.\ J.\ Woldar, \newblock
				%Graphs of prescribed girth and bi-degree, \newblock
				%{\em J. Comb. Theory Ser. B  } {\bf 64} (1995) 228--239.
				
				\bibitem{hl63}
				H. Halberstam and R. Laxton,
				On perfect difference sets, 
				{\em Quart. J. Math.} {\bf 14} (1963), no. 1, 86--90.
				
				%\bibitem{Hur16} M. R.\ Hurley.
				%New geometric large sets.
				%{\em Dissertation thesis. Florida, Atlantic University}  (2016).
				
				%\bibitem{Jaj&Sir}
				%R.\ Jajcay and J.\ \vv  Sir\' a\vv  n.
				%Small vertex-transitive graphs of given degree and girth.
				%{\em Ars Mathematica Contemporanea} {\bf  4} (2011) 375--384.
				
				%\bibitem{Kee}
				%P. Keevash, The existence of designs, arXiv:1401.3665v3 [math.CO] 2 Aug 2019.
				%\bibitem{Sachs63}
				%H. Sachs.
				%Regular graphs with given girth and restricted circuits.
				%{\em J. London Math. Soc. } {\bf 38} (1963) 423--429.
				
				%\bibitem{Tutte47} W. T. Tutte.
				%A family of cubical graphs.
				%{\em  Proc. Cambridge Philos. Soc. } {\bf 43} (1947).
				
				% \bibitem{MP13}
				% B. D. McKay and A. Piperno, Practical Graph Isomorphism, II, 
				% {\em J. Symbolic Computation} {\bf 60} (2013) 94--112.
				
				% \bibitem{VM98}
				% H. van Maldeghem,
				% {\em Generalized Polygons},
				% Birkh\"auser, Springer, Basel (1998).
				
				\bibitem{ms16}
				M. Miller and J. \v{S}ir\'a\v{n},
				Moore graphs and beyond: A survey,
				{\em Electronic J. Combin.\/} \textbf{20 (2)} (2013) \#DS14v21.
				
				\bibitem{js38}
				J. Singer,
				A theorem in finite projective geometry and some applications to number
				theory,
				{\em Trans. Amer. Math. Soc.} {\bf 43} (1938), no. 3, 377--385.
				
				\bibitem{yff83}
				J. L. A. Yebra, M. A. Fiol, and J. F\`abrega,
				Semiregular bipartite Moore graphs,
				{\em Ars Combin.} {\bf 16A} (1983) 131--139.
				
				%\bibitem{Won}
				%P.K.\ Wong.
				%Cages - a survey.
				%{\em J. Graph Theory } {\bf 6} (1982) 1-22.
				
				
			\end{thebibliography}
		\end{document}